\documentclass[11pt]{article}

\usepackage{amsmath}
\usepackage{amssymb}  
\usepackage{epsfig}
\usepackage{amsthm}   
\usepackage[mathcal]{eucal}
\usepackage{url}

\newcommand{\F}{\mathcal{F}}
\newcommand{\J}{\mathcal{J}}
\newcommand{\Sf}{\mathcal{S}}

\newcommand{\R}{\mathbb{R}}
\newcommand{\BR}{\bar{\mathbb{R}}}

\newcommand{\inner}[2]{\langle{#1},{#2}\rangle}

\newcommand{\norm}[1]{\|#1\|}
\newcommand{\Norm}[1]{\left\|#1\right\|}
\newcommand{\normq}[1]{ {\|#1\|}^2 }

\newcommand{\tos}{\rightrightarrows}

\DeclareMathOperator{\clconv}{cl\,conv}
\DeclareMathOperator{\cl}{cl}
\DeclareMathOperator{\conv}{conv}

\newtheorem{theorem}{Theorem}[section]
\newtheorem{lemma}[theorem]{Lemma}
\newtheorem{corollary}[theorem]{Corollary}

\newtheorem{proposition}[theorem]{Proposition}

\newtheorem{definition}{Definition}[section]

\title{On the surjectivity properties of perturbations of maximal
  monotone operators in non-reflexive Banach spaces}

\author{M. Marques Alves\thanks{IMPA, Estrada Dona Castorina 110, 22460-320
    Rio de Janeiro, Brazil
   ({\tt maicon@impa.br})}\hspace{.5em}\thanks{Partially supported by Brazilian CNPq
    scholarship 140525/2005-0.}
  \and
    B. F. Svaiter\thanks{ IMPA, Estrada Dona Castorina 110, 22460-320 Rio de
    Janeiro, Brazil ({\tt benar@impa.br}) }\hspace{.5em}
    \thanks{Partially supported by CNPq
    grants 300755/2005-8, 475647/2006-8 and by PRONEX-Optimization}
}

\date{}

\begin{document}

\maketitle

\begin{abstract}
  We are concerned with surjectivity of perturbations of maximal
  monotone operators in non-reflexive Banach spaces. While in a
  reflexive setting, a classical surjectivity result due to
  Rockafellar gives a necessary and sufficient condition to maximal
  monotonicity, in a non-reflexive space we characterize maximality
  using a ``enlarged'' version of the duality mapping, introduced
  previously by Gossez. 
  \\
\end{abstract}

\pagestyle{plain}

\section{Introduction}

Let $X$ be a real Banach space and $X^*$ its topological dual.  We use the notation $\pi$ and $\pi_*$ for the duality product in $X\times X^*$ and in $X^*\times X^{**}$, respectively:
\begin{align}
 \nonumber
  &\pi:X\times X^*\to\R, &&\pi_*:X^*\times X^{**}\to\R\\
  \label{eq:df.pps}
  &\pi(x,x^*)=\inner{x}{x^*},&  &\pi_*(x^*,x^{**})=\inner{x^*}{x^{**}}.
\end{align}
The norms on  $X$, $X^*$ and $X^{**}$ will be denoted by
$\| \cdot\|$.  We also use the
notation $\BR$ for the extended real numbers:
\[
\BR=\{-\infty\}\cup\R\cup\{\infty\}.
\]
Whenever necessary, we will identify $X$ with its image under the
canonical injection of $X$ into $X^{**}$. 

A point to set operator $T:X\tos X^*$ is a relation on $X\times X^*$:
\[ 
 T\subset X\times X^* 
\]
and $T(x)=\{x^*\in X^*\;|\; (x,x^*)\in T\}$.
An operator $T:X\tos X^*$ is {\it monotone} if
\[
\inner{x-y}{x^*-y^*}\geq 0,\forall (x,x^*),(y,y^*)\in T
\]
and it is {\it maximal monotone}  if it
is monotone and maximal (with respect to the inclusion) in the family
of monotone operators of $X$ into $X^*$.
The conjugate of $f$ is $f^*:X^*\to \BR$,
\[
 f^*(x^*)=\sup_{x\in X} \inner{x}{x^*}-f(x).
\]
Note that $f^*$ is always convex and lower semicontinuous.

The {\it subdifferential} of $f$ is the point to set operator
$\partial f:X \tos X^{*}$ defined at $x\in X$ by
\[ 
 \partial f(x)=\{x^*\in X^{*}\,|\, f(y)\geq f(x)+\inner{y-x}{x^*},\quad \forall y\in X\}.
\]
For each $x\in X$, the elements $x^*\in \partial f(x)$ are called {\it
  subgradients} of $f$.
The concept of {\it $\varepsilon$-subdifferential} of a convex
function $f$ was introduced by Br\o ndsted and Rockafellar
\cite{BroRock65}. It is a point to set operator
$\partial_{\varepsilon}f:X \tos X^{*}$ defined at each $x \in X$ as
\[ 
\partial_{\varepsilon} f(x) = \{ x^*\in X^{*} 
         \,|\, f(y)\geq f(x)+\inner{y-x}{x^*}-
         \varepsilon,
     \quad \forall y\in X \},
\]
where $ \varepsilon \geq 0$. Note that $\partial f = \partial_{0} f$ 
and $\partial f(x) \subset \partial_{\varepsilon} f(x)$, for all $
\varepsilon \geq 0.$ 

A convex function $f:X\to\BR$ is said to be proper if $f>-\infty$ and
there exists a point $\hat x\in X$ for which $f(\hat x)< \infty$. 
Rockafellar proved that if $f$ is proper, convex and lower
semicontinuous, then $\partial f$ is maximal monotone on $X$ \cite{Rock70}.
If $f:X\to \BR$ is proper, convex and lower semicontinuous,
then $f^*$ is proper and $f$ satisfies  {\it Fenchel-Young inequality}:
for all $x\in X$, $x^*\in X^*$,
\begin{equation} \label{eq:F-Y}
  f(x)+f^*(x^*)\geq \inner{x}{x^*},\quad f(x)+f^*(x^*)=\inner{x}{x^*} \iff
  x^*\in \partial f(x).
\end{equation}
Moreover, in this case, 
$\partial_\varepsilon f$ (and $\partial f=\partial_0 f$) may be
characterized using $f^*$:
\begin{equation}
  \label{eq:eps-sub} 
  \begin{array}{rcl}
         \partial f(x)&=&
  \{ x^*\in X^*\,|\,f(x)+f^*(x^*) = \inner{x}{x^*}\},\\[.4em]
     \partial_{\varepsilon} f(x)&=&
  \{ x^*\in X^*\,|\,f(x)+f^*(x^*)\leq \inner{x}{x^*}+\varepsilon\}.
  \end{array}
\end{equation}
The subdifferential and the  $\varepsilon$-subdifferential of the function
$\frac{1}{2}\|\cdot\|^2$  will be of special interest in this paper,
and will be denoted by $J:X\tos X^*$ and $J_\varepsilon:X\tos X^*$ respectively 
\[
J(x)=\partial\; \frac{1}{2}\|x\|^2,\qquad J_\varepsilon
(x)=\partial_\varepsilon\; \frac{1}{2}\|x\|^2.
\]
Using $f(x)=(1/2)\normq{x}$ in~\eqref{eq:eps-sub},
it is trivial to verify that
\begin{eqnarray*}
  J(x)&=&
  \{ x^*\in X^*\,|\,\frac{1}{2}\normq{x}+\frac{1}{2}\normq{x^*} 
  = \inner{x}{x^*}\}\\[.4em]
  &=&\{x^*\in X^*\,|\, \|x\|^2=\|x^*\|^2=\inner{x}{x^*}\}
\end{eqnarray*}
and 
\[ J_{\varepsilon} (x)=
  \{ x^*\in X^*\,|\,\frac{1}{2}\normq{x}+\frac{1}{2}\normq{x^*} 
        \leq \inner{x}{x^*}+\varepsilon\}.
\]
The operator $J$ is widely used in Convex Analysis in Banach spaces
and it is called the {\it duality mapping} of $X$. The operator
$J_{\varepsilon}$ was introduced by Gossez~\cite{GosProc71} to
generalize some results concerning maximal monotonicity in reflexive
Banach spaces to non-reflexive Banach
spaces. It was also used in~\cite{FitzPhelpSet95} to the study of
locally maximal monotone operators in non-reflexive Banach spaces.

If $X$ is a real \emph{reflexive} Banach space and $T:X\tos X^*$ is
monotone, then $T$ is maximal monotone if and only if
\[ R(T(\cdot +z_0)+J)=X^*, \qquad \forall z_0\in X.
\]
We shall prove a similar result for a class of maximal monotone
operators in non-reflexive Banach spaces.

\section{Basic definitions and theory}
\label{sec:pr}

In this section we present the tools and results which will be used to
prove the main results of this paper.

For $f:X\to\BR$, $\conv f:X\to\BR$ is the largest convex function
majorized by $f$, and $\cl f:X\to \BR$ is the largest lower
semicontinuous function  majorized by $f$.  It is trivial to verify that
\[ \cl f(x)=\liminf_{y\to x} f(y),\qquad f^*=(\conv f)^*=(\clconv f)^*.
\]
The functions $\cl f$ and $\clconv f$ are usually called the (lower
semicontinuous) closure of $f$ and the convex lower semicontinuous
closure of $f$, respectively.

Fitzpatrick proved constructively that maximal monotone operators are
representable by  convex functions.
Let  $T:X\tos X^*$ be maximal monotone.  The \emph{Fitzpatrick
  function of $T$}~\cite{Fitz88} is  $\varphi_T: X\times X^*\to \BR$
\begin{equation} \label{FitzIntro}
 \varphi_{T}(x,x^*)=\sup_{(y,
      y^*)\in T} \inner{x - y}{y^* - x^*}+\inner{x}{x^*}
\end{equation}
and \emph{Fitzpatrick family} associated with $T$ is
\begin{equation} \label{eq:def.ft}
  \F_T=\left\{ h\in \BR^{X\times X^*}
    \left|
      \begin{array}{ll}
        h\mbox{ is convex and lower semicontinuous}\\
        \inner{x}{x^*}\leq h(x,x^*),\quad \forall (x,x^*)\in X\times X^*\\
        (x,x^*)\in T 
        \Rightarrow 
        h(x,x^*) = \inner{x}{x^*}
      \end{array}
    \right.
  \right\}.
\end{equation}
\begin{theorem}[\mbox{\cite[Theorem 3.10]{Fitz88}}] \label{th:fitz} 
  Let $X$ be a real Banach space and $T:X\tos X^*$ be maximal
  monotone. Then  for any $h\in \F_T$ \eqref{eq:def.ft}
  \[
  (x,x^*)\in T\iff h(x,x^*)=\inner{x}{x^*}, \qquad \forall (x,x^*)\in
  X\times X^*
  \]
  and   $\varphi_T$ \eqref{FitzIntro} is the smallest element
  of the  family $\F_T$.
\end{theorem}
\noindent
Fitzpatrick's results described above were rediscovered by
Mart\'inez-Legaz and Th\'era~\cite{LegTheJNCA01}, and Burachik and
Svaiter~\cite{BuSvSet02}. Since then, this area has been subject of
intense research.

The \emph{indicator function} of $A\subset X$ is $\delta_A:X\to\BR$,
\[\delta_A(x):=
\begin{cases}
  0,& x\in A\\
  \infty,& \mbox{ otherwise.}
\end{cases}
\]
Using the indicator function we have another expression for
Fitzpatrick function:
\[
 \varphi_T(x,x^*)= \left(\pi+\delta_T\right)^*(x^*,x).
\]
The supremum of Fitzpatrick family is the $\mathcal{S}$-function,
defined and studied by Burachik and Svaiter in~\cite{BuSvSet02},
$\mathcal{S}_T:X\times X^*\to\BR$
\[
\mathcal{S}_T(x,x^*)=\sup  \left\{ h(x,x^*)\;\left|\;
    \begin{array}{l}
  h: X\times X^*\to \BR  \mbox{  convex lower semicontinuous}\\
  h(x,x^*)\leq \inner{x}{x^*}, \quad\forall (x,x^*)\in T
    \end{array}
  \right\}\right.
\]
or, equivalently (see~\cite[Eq.(35)]{BuSvSet02}, \cite[Eq. 29]{BuSvIMPA01})
\begin{equation}
  \label{eq:def.s}
  \mathcal{S}_T=\clconv (\pi+\delta_T).
\end{equation}
Some authors~\cite{BorJCA06,VosSet06,BorProc07} attribute the
$\mathcal{S}$-function to~\cite{PenRelv04} although
this work
was \emph{submitted} after the publication of~\cite{BuSvSet02}.
Moreover, the content of~\cite{BuSvSet02}, and specifically the $\mathcal{S}_T$
function, was presented on Erice workshop
on July 2001, by R. S. Burachik~\cite{BuErice01b}.  A list of the
talks of this congress, which includes~\cite{PenErice01}, is available
on the www\footnote{
  \url{http://www.polyu.edu.hk/~ama/events/conference/EriceItaly-OCA2001/Abstract.html}}.
It shall also be noted that~\cite{BuSvIMPA01}, the preprint
of~\cite{BuSvSet02}, was published ( and available on www) at IMPA
preprint server in August 2001.

Burachik and Svaiter also proved that the family $\F_T$ is invariant under
the mapping 
\begin{equation}
  \label{eq:def.jconj}
  \J:\BR^{X\times X^*}\to\BR^{X\times X^*}, \;\J\;h(x,x^*)=h^*(x^*,x).
\end{equation}
If $T:X\tos X^*$ is maximal monotone, then~\cite{BuSvSet02}
\[ \J(\F_T)\subset \F_T,\quad \J\;\Sf_T=\varphi_T.
\]
In particular, for any $h\in \F_T$,
\begin{equation}
  \label{eq:hhs}
  h(x,x^*)\geq \inner{x}{x^*},\quad h^*(x^*,x)\geq \inner{x}{x^*},
\qquad \forall (x,x^*)\in X\times X^*.
\end{equation}
A partial converse of this fact was proved in~\cite{BuSvProc03}: 
in a \emph{reflexive}
Banach space, if $h$ is convex, lower semicontinuous and 
satisfy 
\eqref{eq:hhs}
then
$$T:=\{(x,x^*)\;|\; h(x,x^*)=\inner{x}{x^*}\}$$
is maximal monotone and $h\in \F_T$~\cite{BuSvProc03}.
In order to extend this result  to non-reflexive Banach
spaces, Marques Alves and Svaiter considered an extension of
condition~\eqref{eq:hhs} to non-reflexive Banach spaces:
\begin{equation}
 \label{eq:hhs.e}
\begin{array}{rll}
  h(x,x^*)&\geq
  \inner{x}{x^*}, &\forall (x,x^*)\in X\times X^*,\\
  h^*(x^*,x^{**})&\geq \inner{x^*}{x^{**}}, & \forall (x^*,x^{**})\in
  X^*\times X^{**}.
\end{array}
\end{equation}
We shall prefer the 
synthetic notation $h\geq \pi$,
$h^*\geq\pi_*$
for the above condition.
The following result will be fundamental in our analysis
\begin{theorem}[\mbox{\cite[Theorem 3.4]{MASvJCA08}}]
  \label{th:br}
  Let $h:X\times X^*\to \BR$ be a convex and lower semicontinuous function. If 
  \[ h\geq \pi,\qquad h^*\geq \pi_*\]
  and $ h(x,x^*)< \inner{x}{x^*}+\varepsilon$, then for any
  $\lambda>0$ there exists $x_\lambda$, $x_\lambda^*$ such that 
  \[ h(x_\lambda,x_\lambda^*)=\inner{x_\lambda}{x_\lambda^*},
    \qquad
    \norm{x_\lambda-x}<\lambda,\quad\norm{x_\lambda^*-x^*}<\varepsilon/\lambda.
  \]
\end{theorem}
Using Theorem~\ref{th:br}, the authors proved~\cite{MASvJCA08} that
condition~\eqref{eq:hhs.e} ensures that $h$ represents a maximal
monotone operator. Here we will be interested also in the case where the
lower semicontinuity assumption is removed.
\begin{theorem}[\mbox{\cite[Theorem 4.2, Corollary 4.4]{MASvJCA08}}]
  \label{th:br2}
  Let $h:X\times X^*\to \BR$ be a convex function. If 
  \[ h\geq \pi,\qquad h^*\geq \pi_*\]
  then 
  \[ T=\{(x,x^*)\in X\times X^*\,|\,h^*(x^*,x)=\inner{x}{x^*}\}
 \]
 is maximal monotone and satisfy the restricted Br\o
 ndsted-Rockafellar property.
 Additionally, if $h$ is also lower semicontinuous, then
 \[ T=\{(x,x^*)\in X\times X^*\,|\,h(x,x^*)=\inner{x}{x^*}\}.
 \]
\end{theorem}
We will need the following immediate consequence of the above theorem:
\begin{corollary}
  \label{cr:cl}
  Let $h:X\times X^*\to \BR$. If 
  \[ \conv h\geq \pi,\qquad h^*\geq \pi_*\]
  then 
  \begin{align*}
     T&=\{(x,x^*)\in X\times X^*\,|\,h^*(x^*,x)=\inner{x}{x^*}\}\\
   &=\{(x,x^*)\in X\times X^*\,|\,\J h(x,x^*)=\inner{x}{x^*}\}
  \end{align*}
  is maximal monotone,
  \[  T=\{(x,x^*)\in X\times X^*\,|\,\clconv h(x,x^*)=\inner{x}{x^*}\}
 \]
  $\clconv h\in \F_T$ and $\J h\in\F_T$, where
  $  \J h(x,x^*)=h^*(x^*,x)$.
\end{corollary}
\begin{proof}
  As the duality product is continuous in $X\times X^*$, $\clconv h\geq
  \pi$. As conjugation is invariant under  the $\conv$ operation
  and the  (lower semicontinuous)
  closure, $(\clconv h)^*=h^*\geq \pi_*$. 
  To end the proof, apply
  Theorem~\ref{th:br2} to $\clconv h$,  observe that $\J h$ is convex,
  lower semicontinuous, $\J h\geq \pi$ and use  definition~\eqref{eq:def.ft}.
\end{proof}

In a non-reflexive Banach Space $X$, if $T:X\tos X^*$ is maximal
monotone and for some $h\in\F_T$ it holds that $h\geq \pi$, $h^*\geq
\pi_*$, then $T$ behaves similarly to a maximal monotone operator in a
\emph{reflexive} Banach space.  A natural question is: what is the
class of maximal monotone operators (in non-reflexive Banach spaces)
which have some function in Fitzpatrick family
satisfying~\eqref{eq:hhs.e}?
To answer this question, first let us recall the definition of maximal monotone operators of type
NI~\cite{SimRanJMA96}.
\begin{definition}
  \label{def:ni}
  A maximal monotone operator $T:X\tos X^*$ is type NI if
  \[ 
  \inf_{(y,y^*)\in T}\inner{y^*-x^*}{x^{**}-y}\leq 0,\qquad
  \forall (x^*,x^{**})\in X^*\times X^{**}.
  \]   
\end{definition}
In~\cite{VoZaArx08} it was observed that if $T$ is a maximal monotone operators
of type NI, then $\Sf_T$ satisfies condition~\eqref{eq:hhs.e}.
We shall need the following theorem. As it is
proved in a paper not yet published, we include its proof on the Appendix~\ref{ap:a}.
\begin{theorem}[\mbox{\cite[Theorem 1.2]{MASvJCA08.2}}]
  \label{th:aux}
  Let $T:X\tos X^*$ be maximal monotone. The following conditions are
  equivalent
  \begin{enumerate}
  \item $T$ is type NI,
  \item there exists $h\in\F_T$ such that $h\geq \pi$ and $h^*\geq \pi_*$,
  \item for all $h\in\F_T$,  $h\geq \pi$ and $h^*\geq \pi_*$,
  \item there exists $h\in\F_T$ such that
    \[
    \inf h_{(x_0,x_0^*)}+\frac{1}{2}\normq{x}+
    \frac{1}{2}\normq{x^*}=0,
    \qquad \forall (x_0,x_0^*)\in X\times X^*,
    \]
  \item for all $h\in\F_T$,
   \[
    \inf h_{(x_0,x_0^*)}+\frac{1}{2}\normq{x}+
    \frac{1}{2}\normq{x^*}=0,
    \qquad \forall (x_0,x_0^*)\in X\times X^*.
    \]
  \end{enumerate}
\end{theorem}

\section{Surjectivity and maximal monotonicity in non-reflexive Banach
  spaces}
\label{sec:sp}

We begin with two elementary technical results which will be useful.	
\begin{proposition}
 \label{prop:tch1}
 The following statements holds:
\begin{enumerate}
 \item For any $\varepsilon\geq 0$, if $y^*\in J_{\varepsilon}(x)$, then 
 $
 \left|\;\Norm{x}-\Norm{y^*}\;\right|\leq \sqrt{2\varepsilon}$.
\item Let $T:X\tos X^*$ be a monotone operator and
  $\varepsilon,M>0$. Then,
\[ \left(T+J_\varepsilon\right)^{-1}(B_{X^*}[0,M]) \]
is bounded.
\end{enumerate}
\end{proposition}
\begin{proof}
  To prove item 1, let $\varepsilon\geq 0$ and $y^*\in
  J_{\varepsilon}(x)$. The desired result follows from the following
  inequalities:
\[
 \frac{1}{2}(\|x\|-\|y^*\|)^2\leq \frac{1}{2}\|x\|^2+\frac{1}{2}\|y^*\|^2-\inner{x}{y^*}
 \leq \varepsilon.
\]
To prove item 2, take $(z,z^*)\in T$. If  
$x\in \left(T+J_\varepsilon\right)^{-1}(B[0,M])
$
then there exists $x^*, y^*$ such that
\[
x^*\in T(x),\quad y^*\in
J_{\varepsilon}(x),\qquad \norm{x^*+y^*}\leq M.
\]
Therefore, using Fenchel Young inequality~\eqref{eq:F-Y}, the monotonicity of $T$ and
the definition of $J_\varepsilon$ we obtain
\begin{align*}
   \frac{1}{2}\normq{x-z}+ \frac{1}{2}\normq{x^*+y^*-z^*}
&\geq \inner{x-z}{x^*+y^*-z^*}\\
&\geq \inner{x-z}{y^*}\\
&\geq\left[
  \frac{1}{2}\normq{x}+\frac{1}{2}\normq{y^*}-\varepsilon\right]
   -\norm{z}\norm{y^*}.
\end{align*}
Note also that
\[ 
\normq{x-z}\leq\normq{x}+2\norm{x}\norm{z}+\normq{z},\qquad
\normq{x^*+y^*-z^*}\leq (M+\norm{z^*})^2.
\]
Combining the above equations we obtain
\begin{align*}
  \frac{1}{2}\normq{z}+\frac{1}{2}(M+\norm{z^*})^2
  &\geq\frac{1}{2}\normq{y^*}-\norm{x}\norm{z}- \norm{z}\norm{y^*}
  -\varepsilon.
\end{align*}
As $y^*\in J_\varepsilon(x)$, by
item 1, we have $\norm{x}\leq\norm{y^*}+\sqrt{2\varepsilon}$.
Therefore
\begin{align*}
  \frac{1}{2}\normq{z}+\frac{1}{2}(M+\norm{z^*})^2
  &\geq\frac{1}{2}\normq{y^*}-2\norm{y^*}\norm{z}-
  \norm{z}\sqrt{2\varepsilon} -\varepsilon.
\end{align*}
Hence, $y^* $ is bounded. In fact,
\[ \norm{y^*}\leq
  2\norm{z}+\sqrt{
4\normq{z}+2\left[
\norm{z}\sqrt{2\varepsilon}+\varepsilon
\right]+\normq{z}+(M+\norm{z^*})^2.
}
\]
As we already observed,
$\norm{x}\leq\norm{y^*}+\sqrt{2\varepsilon}$ and so, $x$ is also
bounded.
\end{proof}

Now we will prove that under monotonicity, dense range of some
perturbation of a monotone operator is equivalent to surjectivity of
that perturbation.
\begin{lemma}
  \label{lm:b1}
  Let $T:X\tos X^*$ be monotone and $\mu>0$. Then the conditions below
  are equivalent
  \begin{enumerate}
  \item  $ \overline{R(T(\cdot+z_0)+\mu J_\varepsilon)}=X^*$,
    for any $\varepsilon>0$ and  $z_0\in X$, 
 \item   $ R(T(\cdot+z_0)+\mu J_\varepsilon)=X^*$
 for any $\varepsilon>0$ and $z_0\in X$. 
  \end{enumerate}
\end{lemma}
\begin{proof}
  It suffices to prove the lemma for $\mu=1$ and then, for the general
  case, consider $T'=\mu^{-1}T$. Now note that for any $z_0\in X$ and
  $z_0^*\in X^*$, $T-\{(z_0,z_0^*)\}$ is also monotone. Therefore, it
  suffices to prove that $0\in \overline{R(T+J_\varepsilon)}$, for any
  $\varepsilon>0$ if and only if $0\in R(T+J_\varepsilon)$, for any
  $\varepsilon>0$.
  The "if" is easy to check. To prove the "only if", suppose that
  \[
  0\in \overline{R(T+J_\varepsilon)},\qquad \forall\varepsilon>0.
  \]
  First use item 2 of Proposition~\ref{prop:tch1} with $M=1/2$ to
  conclude that there exists $\rho>0$ such that
  \[
  (T+J_{1/2})^{-1}\left(B_{X^*}[0,1/2]\right)\subset B_X[0,\rho].
  \]
  By assumption, for any $0<\eta<\frac{1}{2}$ there exists
  $x_{\eta}\in X$, $x_{\eta}^*,y_{\eta}^*\in X^*$ such that
  \begin{equation}
    \label{eq:lm.b1.1}
    x_{\eta}^*\in T(x_{\eta}),\quad y_{\eta}^*\in J_{\eta}(x_{\eta})\quad \mbox{and}
    \quad \|x_{\eta}^*+y_{\eta}^*\|< \eta<\frac{1}{2}.
  \end{equation}
  As $J_\eta(x_\eta)\subset J_{1/2}(x_\eta)$, $x_\eta\in
  (T+J_{1/2})^{-1}(x_{\eta}^*+y_{\eta}^*)$ and so,
  \[ \|x_\eta\|\leq \rho,\qquad \|y_\eta^*\|\leq \rho+1.
  \]
  where the second inequality follows from the first one and item 1 of
  Proposition~\ref{prop:tch1}. Therefore
  \begin{align*}
      \frac{1}{2}\|x_{\eta}^*\|^2&\leq 
    \frac{1}{2}\left(\|x_{\eta}^*+y_{\eta}^*\|+\|y_{\eta}^*\|\right)^2
    \leq \frac{1}{2}\eta^2+\eta(\rho+1)+\frac{1}{2}\|y_{\eta}^*\|^2,\\[.5em]
   \inner{x_\eta}{x^*_\eta}&=\inner{x_\eta}{x_\eta^*+y_\eta^*}-
    \inner{x_\eta}{y_\eta^*}
      \leq\rho\eta- \inner{x_\eta}{y_\eta^*}.
  \end{align*}
  Combining the above inequalities we obtain
  \begin{align*}
    \frac{1}{2}\|x_{\eta}\|^2+\frac{1}{2}\|x_{\eta}^*\|^2+\inner{x_{\eta}}{x_{\eta}^*}
   &\leq  \frac{1}{2}\|x_{\eta}\|^2+\frac{1}{2}\|y_{\eta}^*\|^2
    -\inner{x_{\eta}}{y_{\eta}^*} 
    +\eta(2\rho+1)+\frac{1}{2}\eta^2.
  \end{align*}
  The inclusion
  $y^*_{\eta}\in
  J_{\eta}(x_{\eta})$, means that,
  \begin{equation}
    \label{eq:lm.b1.3}
    \frac{1}{2}\|x_{\eta}\|^2+\frac{1}{2}\|y_{\eta}^*\|^2
     -\inner{x_{\eta}}{y_{\eta}^*}\leq \eta.
  \end{equation}
  Hence, using the two above inequalities we conclude that
\[
 \frac{1}{2}\|x_{\eta}\|^2+\frac{1}{2}\|x_{\eta}^*\|^2+\inner{x_{\eta}}{x_{\eta}^*}
   \leq 2\eta(\rho+1)+\frac{1}{2}\eta^2.
\]
To end the prove, take an arbitrary $\varepsilon>0$. Choosing 
$0<\eta<1/2$ such that,
 \[
  2\eta(\rho+1)+\frac{1}{2}\eta^2<\varepsilon,
\]
we have
\begin{align*}
\frac{1}{2}\|x_{\eta}\|^2+\frac{1}{2}\|x_{\eta}^*\|^2
+\inner{x_{\eta}}{x_{\eta}^*}<\varepsilon,\qquad   x_\eta^*\in T(x_\eta).
\end{align*}
According tho the above inequality,
 $-x_\eta^*\in J_\varepsilon(x_\eta)$. Hence
$0\in(T+J_\varepsilon)(x_\eta)$.
\end{proof}

In a reflexive Banach space, surjectivity of a monotone operator plus
the duality mapping is equivalent to maximal monotonicity. This is a
classical result of Rockafellar~\cite{RockSum70}.
To obtain a
partial extension of this
result to non-reflexive Banach spaces, we must consider the ``enlarged''
duality mapping. 
\begin{lemma}
  \label{lm:b2}
  Let $T:X\tos X^*$ be  monotone and $\mu>0$. If
  \[
  \overline{R(T(\cdot+z_0)+ \mu J_\varepsilon)}=X^*,\qquad 
  \forall \varepsilon>0, z_0\in X\]
  then $\overline T$,  the closure of $T$ in the norm-topology of $X\times
  X^*$, is maximal monotone and type NI.
\end{lemma}
\begin{proof}
  Note that $ T+\mu J_\varepsilon=\mu(\mu^{-1}T+J\varepsilon)$.
  Therefore, it suffices to prove the lemma for $\mu=1$ and then, for
  the general case, consider $T'=\mu^{-1}T$. The monotonicity of $\bar
  T$ follows from the continuity of the duality product.

  Using the assumptions on $T$ and Lemma~\ref{lm:b1} we conclude
  that $T(\cdot+z_0)+J_\varepsilon$ is onto, for any $\varepsilon>0$
  and $z_0\in X$. Therefore, for any $(z_0,z_0^*)\in X\times X^*$ and
  $\varepsilon>0$, there exists $x_\varepsilon$, $x^*_\varepsilon$ such that
  \begin{equation}
    \label{eq:lm.b2.2}
    x_\varepsilon^*+z_0^*\in T(x_\varepsilon+z_0)\quad \text{and } -x_\varepsilon^*\in
    J_\varepsilon (x_\varepsilon).
\end{equation}
Note that the second inclusion in the above equation is equivalent to
\begin{equation}
 \label{eq:lm.b2.22}
 \frac{1}{2}\|x_\varepsilon\|^2+\frac{1}{2}\|x_\varepsilon^*\|^2
 \leq \inner{x_\varepsilon}{-x_\varepsilon^*}+ \varepsilon.
\end{equation}

  To prove maximal monotonicity of $\bar T$, suppose that $(z_0,
  z_0^*)\in X\times X^*$ is monotonically related to $\bar T$. As
  $T\subset \bar T$
  \[
    \inner{z-z_0}{z^*-z_0^*}\geq 0,\; \forall\; (z,z^*)\in T.
  \]
  So, taking $\varepsilon>0$ and $x_\varepsilon\in X$,
  $x_\varepsilon^*\in X^*$ as in \eqref{eq:lm.b2.2} we conclude that
  \[
    \inner{x_\varepsilon}{x_\varepsilon^*}=
    \inner{x_\varepsilon+z_0-z_0}{x^*_\varepsilon+z^*_0-z_0^*}\geq 0,
  \]
which, combined with \eqref{eq:lm.b2.22} yields
\[ \frac{1}{2}\normq{x_\varepsilon}+\frac{1}{2}\normq{x^*_\varepsilon}\leq
\varepsilon.
\] 
As $(x_\varepsilon+z_0,x^*_\varepsilon+z_0^*)\in T$, and $\varepsilon$
is an arbitrary strictly positive number, we conclude that
$(z_0,z_0^*)\in\bar T$, and $\bar T$ is maximal monotone.

It remains to prove that $\bar T$ is type NI.
Consider an arbitrary $(z_0,z_0^*)\in X\times X^*$ and $h\in \F_{\bar
  T}$. Then, using~\eqref{eq:lm.b2.2}, \eqref{eq:lm.b2.22} we conclude
that for any $\varepsilon >0$, there exists $(x_\varepsilon,x_\varepsilon^*)\in X\times X^*$ such
that
\[ 
h(x_{\varepsilon}+z_0,x_{\varepsilon}^*+z_0^*)=\inner{x_\varepsilon+z_0}{x_\varepsilon^*+z_0^*}, 
\qquad\frac{1}{2}\|x_\varepsilon\|^2+\frac{1}{2}\|x_\varepsilon^*\|^2
\leq \inner{x_\varepsilon}{-x_\varepsilon^*}+ \varepsilon.
\]
The first equality above is equivalent to 
$h_{(z_0,z_0^*)}(x_\varepsilon,x_\varepsilon^*)=\inner{x_\varepsilon}{x^*_\varepsilon}$.
Therefore,
\[
 h_{(z_0,z_0^*)}(x_\varepsilon,x_\varepsilon^*)+\frac{1}{2}\|x_\varepsilon\|^2+
 \frac{1}{2}\|x_\varepsilon^*\|^2<
\varepsilon,
\]
that is,
 \[
     \inf h_{(z_0,z_0^*)}(x,x^*)+\frac{1}{2}\normq{x}+
    \frac{1}{2}\normq{x^*}=0.
 \]
Now, use item 5 of Theorem~\ref{th:aux} to conclude that $\bar T$ is type NI.
\end{proof}

Direct application of Lemma~\ref{lm:b2} gives the next corollary.
\begin{corollary}
  \label{cr:b2}
   If $T:X\tos X^*$ is monotone, closed, $\mu>0$ and 
  \[
  \overline{R(T(\cdot+z_0)+ \mu J_\varepsilon)}=X^*,\qquad 
  \forall \varepsilon>0, z_0\in X
  \]
  then $T$,  is maximal monotone and type NI.
\end{corollary}
\begin{proof}
  Use Lemma~\ref{lm:b2} and the assumption $T=\bar T$.
\end{proof}

\begin{lemma}
  \label{lm:ni.ni}
  Let $T_1,T_2:X\tos X^*$ be maximal monotone and type NI. Take
  \[
  h_1\in\F_{T_1},\qquad h_2\in \F_{T_2}
  \]
  and define 
  \begin{align*}
  &h:X\times X^*\to\BR\\   
  &h(x,x^*)= \left(h_1(x,\cdot)\Box h_2(x,\cdot)\right)(x^*)
            =\inf_{y^*\in X^*} h_1(x,y^*)+h_2(x,x^*-y^*),
  \end{align*}
   \[D_X(h_i)=\{x\in X\;|\; \exists\; x^*,\quad
   h_i(x,x^*)<\infty\},\qquad i=1,2.
  \]

  If
  \begin{equation}\label{eq:cond} 
  \bigcup_{\lambda>0}\lambda(D_X(h_1)-D_X(h_2))
  \end{equation}
  is a closed subspace then
  \[ h\geq\pi, h^*\geq \pi_*,\qquad
    \J h\geq \pi,(\J h)^*\geq \pi_*,\]
     \begin{align*}
       T_1+T_2&=\{ (x,x^*)\;|\; \J h(x,x^*)=\inner{x}{x^*}\}\\
    &=\{ (x,x^*)\;|\; h (x,x^*)=\inner{x}{x^*}\}
     \end{align*}
  and $T_1+T_2$ is maximal monotone type NI and
  \[ \J h,\cl h\in\F_{T_1+T_2}.\]
\end{lemma}
\begin{proof}
Since $h_1\in \F_{T_1}$ and $h_2\in \F_{T_2}$, $h_1\geq \pi$ and
$h_2\geq \pi$. So
\[ h_1(x,y^*)+h_2(x,x^*-y^*)\geq \inner{x}{y^*}+\inner{x}{x^*-y^*}=
   \inner{x}{x^*}.
\]
Taking the $\inf$ in $y^*$ at the left-hand side of the above
inequality we conclude that $h\geq\pi$.

Let $(x^*,x^{**})\in X^*\times X^{**}$. Using the definition of $h$ we
have
\begin{align}
  h^*(x^*,x^{**})&=
  \sup_{(z,z^*)\in X\times X^*}\inner{z}{x^*}+\inner{z^*}{x^{**}}-h(z,z^*)\\[.5em]
  &= \sup_ {(z,z^*,y^*)\in X\times X^*\times X^*}
  \begin{array}[t]{c}
     \inner{z}{x^*}+\inner{z^*}{x^{**}}-h_1(z,y^*)\\-h_2(z,z^*-y^*) 
  \end{array}\\
  &= \sup_ {(z,y^*,w^*)\in X\times X^*\times X^*}
  \begin{array}[t]{c}
     \inner{z}{x^*}+\inner{y^*}{x^{**}}+\inner{w^*}{x^{**}}-h_1(z,y^*)\\-h_2(z,w^*) 
  \end{array}
\end{align}
where we used the substitution $z^*=w^*+y^*$ in the last term.
So, defining  $H_1,H_2:X\times X^{*}\times X^{*}\to \BR$
\begin{equation}
  \label{eq:def.h12}
  H_1(x,y^*,z^*)=h_1(x,y^*),\quad H_2(x,y^*,z^*)=h_2(x,z^*). 
\end{equation}
we have 
\[ h^*(x^*,x^{**})= (H_1+H_2)^*(x^*,x^{**},x^{**}).
\]
Using~\eqref{eq:cond}, the Attouch-Brezis extension~\cite[Theorem 1.1]{attouch-brezis} of
Fenchel-Rockafellar duality theorem and~\eqref{eq:def.h12} we
conclude that the conjugate of the sum at the right hand side of the
above equation is the \emph{exact} inf-convolution of the
conjugates. Therefore,
\[ h^*(x^*,x^{**})=\min_{(u^*,y^{**},z^{**})}H_1^*(u^*,y^{**},z^{**})
       +H_2^*(x^*-u^*,x^{**}-y^{**},x^{**}-z^{**}).
\]
Direct use of definition~\eqref{eq:def.h12} yields
\begin{equation}
 \label{eq:ni.1}
  H_1^*(u^*,y^{**},z^{**})= h_1^*(u^*,y^{**})+\delta_{0}(z^{**}),\quad
 \forall (u^*,y^{**},z^{**})\in X^*\times X^{**}\times X^{**},
\end{equation}
\begin{equation}
 \label{eq:ni.2}
  H_2^*(u^*,y^{**},z^{**})= h_2^*(u^*,z^{**})+\delta_{0}(y^{**}),\quad
  \forall (u^*,y^{**},z^{**})\in X^*\times X^{**}\times X^{**}.
\end{equation}
Hence,
\begin{equation}
  \label{eq:ni.4}
    h^*(x^*,x^{**})=\min_{u^*\in X^*}h_1^*(u^*,x^{**})
  +h_2^*(x^*-u^*,x^{**}).
\end{equation}
Therefore, using that $h_1^*\geq \pi_*$, $h_2^*\geq
\pi_*$,~\eqref{eq:ni.4} and the same reasoning used to show that
$h\geq \pi$ we have
\[
 h^*\geq \pi^*.
\]
Up to now, we proved that $h\geq \pi$ and $h^*\geq \pi_*$( and $\J h\geq \pi$). So, using Theorem~\ref{th:br2} we conclude that $S:X\tos X^*$, defined as
\[
 S=\{(x,x^*)\in X\times X^*\,|\,\J h(x,x^*)=\inner{x}{x^*}\},
\]
is maximal monotone. As $\J h$ is convex and lower semicontinuous, $\J h\in\F_S$.

 We will prove that $T_1+T_2=S$.  Take $(x,x^*)\in
S$, that is, $\J h(x,x^*)=\inner{x}{x^*}$. Using~\eqref{eq:ni.4} we
conclude that there exists $u^*\in X^*$ such that
\begin{align*}
  h_1^*(u^*,x)+h_2^*(x^*-u^*,x)=\inner{x}{x^*}.
\end{align*}
We know that 
\[ h^*_1(u^*,x)\geq \inner{x}{u^*},
\qquad h_2^*(x^*-u^*,x)\geq \inner{x}{x^*-u^*}.
\]
Combining these inequalities with the previous equation we conclude
that these inequalities holds as equalities, and so
\begin{align*}
  u^*&\in T_1(x),& x^*-u^*\in T_2(x),&& x^*\in (T_1+T_2)(x).\\
 h_1(x,u^*)&=\inner{x}{u^*},
 &h_2(x,x^*-u^*)=\inner{x}{x^*-u^*}, 
 && h(x,x^*)\leq \inner{x}{x^*}. 
\end{align*}
We proved that $S\subset T_1+T_2$. Since $T_1+T_2$ is monotone
and $S$ is maximal monotone, we have $T_1+T_2=S$ (and
$\J h\in\F_{T_1+T_2}$).  Note also that $h(x,x^*)\leq \inner{x}{x^*}$ for any
$(x,x^*)\in T_1+T_2=S$. As $h\geq \pi$, we have equality in $T_1+T_2$.
Therefore,
\[ T_1+T_2\subset \{(x,x^*)\;|\; h(x,x^*)=\inner{x}{x^*}\}
    \subset \{(x,x^*)\;|\;\cl h(x,x^*)\leq \inner{x}{x^*}\}.
\]
Since $h\geq\pi$ and the duality product $\pi$ is \emph{continuous} in
$X\times X^*$, we also have $\cl h\geq\pi$. Hence, using the
above inclusion we conclude that $\cl h$ coincides with $\pi$ in
$T_1+T_2$. Therefore, $\cl h\in\F_{T_1+T_2}$ and the rightmost set in
the above inclusions is $T_1+T_2$. Hence
\[ T_1+T_2=\{(x,x^*)\;|\; h(x,x^*)=\inner{x}{x^*}\}.\]

Conjugation is invariant under the (lower semicontinuous) closure
operation. Therefore,
\[ (\cl h)^*=h^*\geq \pi_*\] and so $T_1+T_2$ is NI.  We proved already
that $\J h\in \F_{T_1+T_2}$. Using item 3 of Theorem~\ref{th:aux} we
conclude that $(\J h)^*\geq \pi_*$.

\end{proof}

\begin{theorem}
  \label{th:1}
  If $T:X\tos X^*$ is a closed monotone operator then the conditions
  bellow are equivalent
  \begin{enumerate}
  \item  $ \overline{R(T(\cdot+z_0)+ J)}=X^*$ for all 
    $z_0\in X$,
  \item $ \overline{R(T(\cdot+z_0)+ J_\varepsilon)}=X^*$ for all 
    $\varepsilon>0$, $z_0\in X$,
  \item $R(T(\cdot+z_0)+ J_\varepsilon)=X^*$ for all 
    $\varepsilon>0$, $z_0\in X$,
  \item $T$ is maximal monotone and type NI.
  \end{enumerate}
\end{theorem}
\begin{proof}
  Item 1 trivially implies item 2. Using Lemma~\ref{lm:b1} we conclude
  that, in particular, item 2 implies item 3. Now use
  Corollary~\ref{cr:b2} to conclude that item 3 implies item 4.  Up to
  now we have 1$\Rightarrow$2$\Rightarrow$3$\Rightarrow$4. 

  To complete the proof we will show that item 4 implies item 1.  So,
  assume that item 4 holds, that is, $T$ is type NI.  
   Take $z_0^*\in X^*$ and $z_0\in X$. 
  Define $T_0=T-\{(z_0,z_0^*)\}$. Trivially
  \[ z_0^*\in \overline{R(T(\cdot+z_0)+ J)}\iff
    0\in \overline{R(T_0+ J)}.
\]
 As the class NI is invariant under translations, 
 in order to prove item 1, it is sufficient to
  prove that if $T$ is type NI, then $0\in \overline{R(T+J)}$. Let
  $h\in \F_{T}$ and $\varepsilon>0$.
  Define $p:X\times X^*\to \R$, 
  \begin{equation} \label{eq:th1.c}
    p(x,x^*)=\frac{1}{2}\|x\|^2+\frac{1}{2}\|x^*\|^2.
  \end{equation}
  Item 5 of
  Theorem~\ref{th:aux} ensure us that there exists
  $(x_{\varepsilon},x_{\varepsilon}^*)\in X\times X^*$ such that
  \begin{equation}
    \label{eq:th1.a}
    h(x_{\varepsilon},x_{\varepsilon}^*)+p(x_{\varepsilon},-x_{\varepsilon}^*)< \varepsilon^2.
  \end{equation}
  Direct calculations yields $p\geq \pi$ and $p^*\geq \pi_*$. We also
  know that $p\in \F_J$ and so $J$ is type NI. Define $H:X\times X^*\to \BR$,
  \[
  H(x,x^*)=\inf_{y^*\in X^*}h(x,y^*)+p(x,x^*-y^*).
  \]
  As $D(p)=X\times X^*$, we
  may apply Lemma~\ref{lm:ni.ni} to conclude that $T+J$ is NI and
  $\cl H\in \F_{T+J}$.
  Using~\eqref{eq:th1.a} we have
  \[
  H(x_\varepsilon,0)\leq h(x_\varepsilon, x_\varepsilon^*)+
  p(x_\varepsilon,-x_\varepsilon^*)<\varepsilon^2.
  \]
  So, $ \cl H(x_\varepsilon,0)\leq H(x_\varepsilon,0)<
  \inner{x_\varepsilon}{0}+\varepsilon^2$. Now use Theorem~\ref{th:br} to
  conclude that there exists $\bar x$, $\bar x^*$ such that
  \[ (\bar x,\bar x^*)\in T+J, \qquad \norm{\bar
    x-x_\varepsilon}<\varepsilon,\qquad 
    \norm{\bar x^*-0}< \varepsilon.
    \]
 So, $\bar x^*\in R(T+J)$ and $\norm{\bar x^*}<\varepsilon$. As
 $\varepsilon>0$ is arbitrary, $0$ is in the closure of $R(T+J)$.
\end{proof}

\begin{corollary}
   If $T:X\tos X^*$ is a closed monotone operator then the conditions
  bellow are equivalent
  \begin{description}
  \item[a]  $ \overline{R(T(\cdot+z_0)+ \mu J)}=X^*$ for all 
    $z_0\in X$ and some $\mu>0$,
  \item[b]  $ \overline{R(T(\cdot+z_0)+ \mu J)}=X^*$ for all 
    $z_0\in X$, $\mu>0$,
  \item[c] $ \overline{R(T(\cdot+z_0)+ \mu J_\varepsilon)}=X^*$ for all 
    $\varepsilon>0$, $z_0\in X$ and some $\mu>0$,
  \item[d] $ \overline{R(T(\cdot+z_0)+ \mu J_\varepsilon)}=X^*$ for all 
    $\varepsilon>0$, $z_0\in X$, $\mu>0$,
 \item[e] $ R(T(\cdot+z_0)+ \mu J_\varepsilon)=X^*$ for all 
    $\varepsilon>0$, $z_0\in X$, and some $\mu>0$,
 \item[f] $ R(T(\cdot+z_0)+ \mu J_\varepsilon)=X^*$ for all 
    $\varepsilon>0$, $z_0\in X$,  $\mu>0$,
  \item [g] $T$ is maximal monotone and type NI.
  \end{description}
\end{corollary}
\begin{proof}
  Suppose that item {\bf a} holds. Define $T'=\mu^{-1}T$ and use
  Theorem~\ref{th:1} to conclude that $T'$ is maximal monotone and
  type NI. Therefore, $T=\mu T'$ is maximal monotone and type NI,
  which means that {\bf g} holds.

  Now assume that item {\bf g} holds, that is, $T$ is maximal monotone
  and type NI. Then, for all $\mu>0$, $\mu^{-1}T$ is maximal monotone
  and type NI, which implies item {\bf b}.

  As the implication {\bf b}$\Rightarrow${\bf a} is trivial, we
  conclude that items {\bf a}, {\bf b}, {\bf g} are equivalent.

  The same reasoning shows that items {\bf c}, {\bf d}, {\bf g} are
  equivalent and so on.
\end{proof}
\appendix

\section{Proof of Theorem~\ref{th:aux}}
\label{ap:a}
 In~\cite{LegSvSet05} Mart\'inez-Legaz and Svaiter defined (with a different notation), 
 for $h:X\times X^*\to\BR$ and
 $(x_0,x_0^*)\in X\times X^*$
 \begin{equation}
  \label{eq:def.hx}
  \begin{array}{l}
    h_{(x_0,x_0^*)}:X\times X^*\to\BR,\\[.4em]
    h_{(x_0,x_0^*)}(x,x^*):=h(x+x_0,x^*+x_0^*)-[\inner{x}{x_0^*}+\inner{x_0}{x^*}
    +\inner{x_0}{x_0^*}].
  \end{array}
 \end{equation}
 The operation $h\mapsto h_{(x_0,x_0^*)}$ preserves many properties of
 $h$, as convexity, lower semicontinuity and can be seen as the action
 of the group $(X\times X^*,+)$ on $\BR^{X\times X^*}$, because
 \[ \left(h_{(x_0,x_0^*)}\right)_{(x_1,x_1^*)}=h_{(x_0+x_1,x_0^*+x_1^*)}.\]
 Moreover
 \[
  \left(h_{(x_0,x_0^*)}\right)^*=\left(h^*\right)_{(x_0^*,x_0)},
 \]
 where the rightmost $x_0$ is identified with its image under the
 canonical injection of $X$ into $X^{**}$.
 Therefore, 
 \begin{enumerate}
 \item $h\geq \pi\iff h_{(x_0,x_0)}\geq \pi$,
 \item $
  \left(h_{(x_0,x_0^*)}\right)^*\geq \pi_*\iff
  \left(h^*\right)_{(x_0^*,x_0)}
  \geq \pi_*$, 
 \end{enumerate}
 The proof of Theorem~\ref{th:aux} will be heavily based on these nice properties of the 
 map $h\mapsto h_{(x_0,x_0^*)}$.

 \begin{proof}[Proof of Theorem~\ref{th:aux}]
 First let us prove that item 2 and item 4 are equivalent. So, suppose item 2 holds and 
 let  $(x_0,x_0^*)\in X\times X^*$. Direct calculations yields
 \[
  h_{(x_0,x_0^*)}\geq \pi,\quad (h_{(x_0,x_0^*)})^*\geq \pi_*.
 \]
 Using~\cite[Theorem 3.1, eq. (12)]{MASvJCA08} we conclude that
 condition item 4 holds.
 For proving that item 4$\Rightarrow$item 2, first note that, for
 any $(z,z^*)\in X\times X^*$,
 \[ h_{(z,z^*)}(0,0)\geq\inf_{(x,x^*)}
 h_{(z,z^*)}(x,x^*)+\frac{1}{2}\normq{x}+ \frac{1}{2}\normq{x^*}.
 \]
 Therefore, using item 4 we obtain
 \[
 h(z,z^*)-\inner{z}{z^*}= h_{(z,z^*)}(0,0)\geq 0.
 \]
 Since $(z,z^*)$ is an arbitrary element of $X\times X^*$ we
 conclude that $h\geq \pi$.
 
 For proving that, $h^*\geq\pi_*$,
 take some $(y^*,y^{**})\in X^*\times X^{**}$.
 First, use Fenchel-Young inequality to conclude that
 for any $(x,x^*), (z,z^*)\in X\times X^*$,
   \begin{align*}
     h_{(z,z^*)}(x,x^*)
     \geq&
     \inner{x}{y^*-z^*}+\inner{x^*}{y^{**}-z}-\left(h_{(z,z^*)}\right)^*(y^*-z^*,y^{**}-z).
   \end{align*}
   As $\left(h_{(z,z^*)}\right)^*=(h^*)_{(z^*,z)}$,
   \begin{align*}
     \left(h_{(z,z^*)}\right)^*(y^*-z^*,y^{**}-z)&=
     h^*(y^*,y^{**})-\inner{z}{y^*-z^*}-\inner{z^*}{y^{**}-z}-\inner{z}{z^*}\\
     &=h^*(y^*,y^{**})-\inner{y^*}{y^{**}}+\inner{y^*-z^*}{y^{**}-z}.
   \end{align*}
   Combining the two above equations we obtain
   \begin{align*}
      h_{(z,z^*)}(x,x^*)
     \geq&
     \inner{x}{y^*-z^*}+\inner{x^*}{y^{**}-z}\\
      &-\inner{y^*-z^*}{y^{**}-z}+\inner{y^*}{y^{**}}-h^*(y^*,y^{**}).
   \end{align*}
   Adding $(1/2)\normq{x}+(1/2)\normq{x^*}$ in both sides of the
   above inequality we have 
  \begin{align*}
      h_{(z,z^*)}(x,x^*)+\frac{1}{2}\normq{x}+\frac{1}{2}\normq{x^*}
     \geq&
     \inner{x}{y^*-z^*}+\inner{x^*}{y^{**}-z}
     +\frac{1}{2}\normq{x}+\frac{1}{2}\normq{x^*}
     \\
      &-\inner{y^*-z^*}{y^{**}-z}+\inner{y^*}{y^{**}}-h^*(y^*,y^{**}).
   \end{align*}
   Note that
   \[
   \inner{x}{y^*-z^*}+\frac{1}{2}\normq{x}\geq
   -\frac{1}{2}\normq{y^*-z^*},\qquad
   \inner{x^*}{y^{**}-z}+\frac{1}{2}\normq{x^*}\geq
   -\frac{1}{2}   \normq{y^{**}-z}.
   \]
   Therefore, for any $(x,x^*), (z,z^*)\in X\times X^*$,
   \begin{align*}
     h_{(z,z^*)}(x,x^*)+\frac{1}{2}\normq{x}+\frac{1}{2}\normq{x^*}
     \geq&
    -\frac{1}{2}\normq{y^*-z^*} -\frac{1}{2}   \normq{y^{**}-z}
     \\
     &-\inner{y^*-z^*}{y^{**}-z}+\inner{y^*}{y^{**}}-h^*(y^*,y^{**}).
   \end{align*}
   Using now the assumption  we conclude that the infimum,
   for $(x,x^*)\in X\times X^*$, at the left hand side of the above
   inequality is $0$. Therefore, taking the infimum on $(x,x^*)\in
   X\times X^*$ at the left hand side of the above inequality and
   rearranging the resulting inequality we have
   \begin{align*}
     h^*(y^*,y^{**})-\inner{y^*}{y^{**}}\geq 
     -\frac{1}{2}\normq{y^*-z^*} -\frac{1}{2}   \normq{y^{**}-z}
     -\inner{y^*-z^*}{y^{**}-z}.
   \end{align*}
   Note that
   \[ \sup_{z^*\in X^*}  -\inner{y^*-z^*}{y^{**}-z}
   -\frac{1}{2}\normq{y^*-z^*}
   =\frac{1}{2}\normq{y^{**}-z}.
  \]
  Hence, taking the sup in $z^*\in X^*$ at the  right hand
  side of the previous inequality we obtain
  \[  h^*(y^*,y^{**})-\inner{y^*}{y^{**}}\geq 0\]
  and item 4 holds.
  Now, using that item 2 and item 4 are equivalent it is trivial to verify that item 3 and item 5 
  are equivalent.

  The second step is to prove that item 4 and item 5 are equivalent.
  So, assume that item 4 holds, that is,
  for some $h\in \F_T$,
  \[
  \inf_{(x,x^*)\in X\times X^*} h_{(x_0,x_0^*)}(x,x^*)+\frac{1}{2}\normq{x}+
  \frac{1}{2}
  \normq{x^*}=0,\qquad \forall (x_0,x_0^*)\in X\times X^*.
  \]
  Take $g\in \F_T$, and $(x_0,x_0^*)\in X\times X^*$.
  First observe that, for any $(x,x^*)\in X\times X^*$, $
  g_{(x_0,x_0^*)}(x,x^*)\geq\inner{x}{x^*}$ and
  \[  g_{(x_0,x_0^*)}(x,x^*)+\frac{1}{2}\normq{x}+
  \frac{1}{2}\normq{x^*}\geq
   \inner{x}{x^*}+\frac{1}{2}\normq{x}+
   \frac{1}{2}\normq{x^*}\geq 0.
   \]
   Therefore,
   \begin{equation}
     \label{eq:aux0}
      \inf_{(x,x^*)\in X\times X^*} g_{(x_0,x_0^*)}(x,x^*)+\frac{1}{2}\normq{x}+
  \frac{1}{2}
  \normq{x^*}\geq 0. 
   \end{equation}
  As the square of the norm is coercive, there exist $M>0$ such that
  \[
  \left\{ (x,x^*)\in X\times X^*\;|\;  h_{(x_0,x_0^*)}(x,x^*)+\frac{1}{2}\normq{x}
    +  \frac{1}{2}\normq{x^*}<1
    \right\}\subset
   B_{X\times X^*}(0,M),
  \]
  where
  \[
  B_{X\times X^*}(0,M)=\left\{ (x,x^*)\in X\times X^*\;|\;
    \sqrt{\normq{x}+\normq{x^*}}<M\right\}.
  \]
  For any $\varepsilon>0$, there exists $(\tilde x,\tilde x^*)$ such
  that
  \[ \min\left\{1,\varepsilon^2\right\}>
  h_{(x_0,x_0^*)}(\tilde x,\tilde x^*)+\frac{1}{2}\normq{\tilde x}
  +  \frac{1}{2}\normq{\tilde x^*}.
  \]
  Therefore
  \begin{equation}
    \label{eq:aux1}
    \begin{array}{l}
      \varepsilon^2 >
      h_{(x_0,x_0^*)}(\tilde x,\tilde x^*)+\frac{1}{2}\normq{\tilde x}
      +\frac{1}{2}\normq{\tilde x^*}\geq
      h_{(x_0,x_0^*)}(\tilde x,\tilde x^*)
       -\inner{\tilde x}{\tilde x^*}\geq 0,\\[.5em]
     M^2\geq \normq{\tilde x}+\normq{\tilde x^*}.
    \end{array}
  \end{equation}
  In particular,
  \[  \varepsilon^2 > h_{(x_0,x_0^*)}(\tilde x,\tilde x^*)
  -\inner{\tilde x}{\tilde x^*}.
  \]
  Now using Theorem~\ref{th:br} we conclude that there exists
  $(\bar x,\bar x^*)$ such that
  \begin{equation}
    \label{eq:aux2}
   h_{(x_0,x_0^*)}(\bar x,\bar x^*)=\inner{\bar x}{\bar x^*},\quad
  \norm{\tilde x-\bar x}<\varepsilon, \quad
  \norm{\tilde x^*-\bar x^*}<\varepsilon.   
  \end{equation}
    Therefore,
   \[ h(\bar x+x_0,\bar x^*+x_0^*)-\inner{\bar x+x_0}{\bar x^*+x_0^*}=
   h_{(x_0,x_0^*)}(\bar x,\bar x^*)-\inner{\bar x}{\bar x^*}=0,
   \]
  and $(\bar x+x_0,\bar x^*+x_0^*)\in T$. As $g\in \F_T$,
  \[
  g(\bar x+x_0,\bar x^*+x_0^*)=\inner{\bar x+x_0}{\bar x^*+x_0^*},
  \]
  and 
  \begin{equation}
    \label{eq:aux3}
    g_{(x_0,x_0^*)}(\bar x,\bar x^*)=\inner{\bar x}{\bar x^*}.
  \end{equation}
  Using the first line of \eqref{eq:aux1} we have
  \[       \varepsilon^2 >
      h_{(x_0,x_0^*)}(\tilde x,\tilde x^*)+
      \bigg[\frac{1}{2}\normq{\tilde x}
      +\frac{1}{2}\normq{\tilde x^*}+
          \inner{\tilde x}{\tilde x^*}
      \bigg]
       -\inner{\tilde x}{\tilde x^*}\geq
\frac{1}{2}\normq{\tilde x}
      +\frac{1}{2}\normq{\tilde x^*}+
          \inner{\tilde x}{\tilde x^*}.
       \]
  Therefore,
  \begin{equation}
    \label{eq:aux4}
     \varepsilon^2>  \frac{1}{2}\normq{\tilde x}
      +\frac{1}{2}\normq{\tilde x^*}+
      \inner{\tilde x}{\tilde x^*}.
  \end{equation}
  Direct use of \eqref{eq:aux2} gives
  \begin{align*}
    \inner{\bar x}{\bar x^*}&=\inner{\tilde x}{\tilde x^*}
    +\inner{\bar x-\tilde x}{\tilde x^*}
    +\inner{\tilde x}{\bar x^*-\tilde x^*}
    +\inner{\bar x-\tilde x}{\bar x^*-\tilde x^*}\\
    &\leq \inner{\tilde x}{\tilde x^*}
    +\norm{\bar x-\tilde x}\,\norm{\tilde x^*}
    +\norm{\tilde x}\,\norm{\bar x^*-\tilde x^*}
    +\norm{\bar x-\tilde x}\,\norm{\bar x^*-\tilde x^*}\\
  &\leq \inner{\tilde x}{\tilde x^*}
    +\varepsilon [ \norm{\tilde x^*}+\norm{\tilde x}]
    +\varepsilon^2
  \end{align*}
  and
  \begin{align*}
    \normq{\bar x}+\normq{\bar x^*}&\leq
    \left( \norm{\tilde x}+\norm{\bar x-\tilde x}\right)^2
    +  \left( \norm{\tilde x^*}+\norm{\bar x^*-\tilde x^*}\right)^2\\
    &\leq \normq{\tilde x}+ \normq{\tilde x^*}
    +2\varepsilon[\norm{\tilde x}+\norm{\tilde x^*}]+2\varepsilon^2
  \end{align*}
  Combining the two above equations with~\eqref{eq:aux3} we obtain
  \[
  g_{(x_0,x_0^*)}(\bar x,\bar x^*)
  +\frac{1}{2}\normq{\bar x}
  +\frac{1}{2}\normq{\bar x^*}\leq
  \inner{\tilde x}{\tilde x^*}
  +\frac{1}{2}\normq{\tilde  x}
  +\frac{1}{2}\normq{\tilde  x^*}+2\varepsilon[\norm{\tilde x}+\norm{\tilde x^*}]+2\varepsilon^2
  \]
  Using now~\eqref{eq:aux4} and the second line of \eqref{eq:aux1} we
  conclude that
 \[
  g_{(x_0,x_0^*)}(\bar x,\bar x^*)
  +\frac{1}{2}\normq{\bar x}
  +\frac{1}{2}\normq{\bar x^*}\leq
  2\varepsilon\;M\sqrt{2}+3\varepsilon^2.
  \]
  As $\varepsilon$ is an arbitrary strictly positive number, using
  also~\eqref{eq:aux0} we conclude that
  \[
      \inf_{(x,x^*)\in X\times X^*} g_{(x_0,x_0^*)}(x,x^*)+\frac{1}{2}\normq{x}+
  \frac{1}{2}
  \normq{x^*}=0.
  \]
  Altogether, we conclude that if item 4 holds then item 5 holds.
  The converse item 5$\Rightarrow$ item 4 is trivial to verify. Hence
  item 4 and item 5 are equivalent. As item 2 is equivalent to item 4 and
  item 3 is equivalent to 5, we conclude that items 2,3,4  and 5 are
  equivalent.

  Now we will prove that item 1 is equivalent to item 3 and conclude the proof
  of the theorem. First suppose that item 3 holds. Since
  $\mathcal{S}_T\in \F_T$
 \[
 (\mathcal{S}_T)^*\geq \pi_*.
 \]
 As has already been observed, for any proper function $h$ it holds 
 that  $(\mbox{cl\,conv}\,h)^*=h^*$.
 Therefore
 \[
 (\mathcal{S}_T)^*=(\pi+\delta_T)^*\geq \pi_*,
 \]
 that is,
 \begin{equation}\label{eq20}
 \sup_{(y,y^*)\in T}
 \inner{y}{x^*}+\inner{y^*}{x^{**}}-\inner{y}{y^*}\geq
 \inner{x^*}{x^{**}}, \forall
    (x^*,x^{**})\in X^*\times X^{**}
 \end{equation}
 After some algebraic manipulations we conclude
 that~\eqref{eq20} is equivalent to
 \[
    \inf_{(y,y^*)\in T}\inner{x^{**}-y}{x^*-y^*}\leq 0,\qquad \forall
    (x^*,x^{**})\in X^*\times X^{**},
 \]
 that is, $T$ is type (NI) and so item 1 holds.
 If item 1 holds, by the same reasoning we conclude that~\eqref{eq20} holds and therefore
 $(\mathcal{S}_T)^*\geq \pi_*$. As $\mathcal{S}_T\in \F_T$, we conclude that item 2 holds.
 As has been proved previously item 2 $\Rightarrow$ item 3.

\end{proof}

\bibliographystyle{plain}

\end{document}